\documentclass{article}

\usepackage{amsmath}
\usepackage[latin1]{inputenc}
\usepackage[T1]{fontenc}
\usepackage[english]{babel}
\usepackage{enumerate}
\usepackage{amssymb}
\usepackage{cleveref}
\usepackage{url}
\usepackage{amsthm}
\usepackage{titlesec} 

\newtheorem{Lemma}{Lemma}
\newtheorem{Theorem}[Lemma]{Theorem}
\newtheorem{Corollary}[Lemma]{Corollary}
\newtheorem{Remark}[Lemma]{Remark}
\newtheorem{Problem}{Problem}
\title{The SOT-closure of the set of 2-isometries}
\author{Marcel Scherer}
\date{}

\begin{document}

\maketitle
\begin{abstract}
We show that the set of $2$-isometries on an infinite-dimensional Hilbert space is not closed in the strong operator topology. In fact, we prove that its SOT-closure coincides with the set of all expansive operators.
\end{abstract}

\let\thefootnote\relax\footnotetext{ \date\\
2020 \textit{Mathematics Subject Classification.} Primary: 47A05, Secondary: 47A20.\\
\ \textit{Keywords and phrases:} strong operator topology; SOT-closure; $2$-isometry; $m$-isometry; Brownian unitary; expansive operator}

Let $H$ be an infinite-dimensional Hilbert space. A bounded linear operator $T\in\mathcal{B}(H)$ is called a \textit{2-isometry} if 
  \begin{equation*}
    1-2T^*T+{T^*}^2T^2=0,
   \end{equation*}
see \cite{JiAg}. If $(T_n)_n$ is a SOT-convergent sequence of $2$-isometries on $H$ with limit $T$, then $(\|T_n\|)_n)$ is bounded by the Banach-Steinhaus theorem and therefore $(T_n^2)_n$ converges SOT to $T^2$. Thus
   \begin{equation*}
     \langle (1-2T^*T+{T^*}^2T^2)x,x\rangle=\lim_{n\to\infty} \|x\|^2-2\|T_nx\|^2+\|T_n^2x\|^2=0,
  \end{equation*}
showing that $T$ is again a $2$-isometry, see \cite{JaZe}. In the same work, the following problem was posed:
  \begin{Problem}
  Is the set of $2$-isometries on an infinite-dimensional Hilbert space strongly closed?
  \end{Problem}
Note that this is not answered by the above observation, because SOT-convergent nets need not be pointwise bounded. A classical example of a set of operators that is SOT-closed under sequences but not under nets is due to Halmos:
  \begin{equation*}
    \{T\in\mathcal{B}(H);\ T^2=0\}.
  \end{equation*}
The SOT-closure of this set equals in fact $\mathcal{B}(H)$, see, for example, \cite{PaHa}.\\
In the first theorem we answer the above problem in the negative. The proof constructs \textit{Brownian unitaries}, i.e. operators of the form
  \begin{equation*}
    \begin{pmatrix} R & \sigma V\\ 0 & U\end{pmatrix},
  \end{equation*}
for $\sigma\in (0,\infty)$, and a decomposition $H_1\oplus H_2=H$ such that $R$ and $V$ are isometries, $\textup{Im}(V)\oplus \textup{Im}(R)=H_1$ and $U$ is unitary. It is easy to check that Brownian unitaries are $2$-isometries.\\
The operator that lies in the closure but is not itself a $2$-isometry will be $2id$. Afterwards we show that every \textit{expansive operator} (i.e. $T$ with $\|Tx\|\ge\|x\|$ for all $x$) is an SOT-limit of $2$-isometries. We do not need the former result for this second statement, however, the main idea behind the proof is most transparent for the operator $2id$.

\begin{Theorem}\label{Theorem1}
Let $H$ be an infinite-dimensional Hilbert space. It holds that
  \begin{equation*}
    2id\in\overline{\{T\in\mathcal{B}(H);\ T\ \textup{2-isometry}\}}^{SOT}.
  \end{equation*}
\end{Theorem}

\underline{Proof}:\\
Let $F\subset H$ be a finite-dimensional subspace. Then there exists $n\in\mathbb{N}$ and an ONB $x_1,\dots,x_n$ of $F$. Set $\epsilon=\tfrac{1}{n}$. Note that $n$, and hence $\epsilon$, depend only on the $\dim(F)$.\\
We claim that there exists an ONS $y_1^{(1)}, y_1^{(2)}, y_2^{(1)},\dots, y_n^{(2)}$ in $H$ such that
  \begin{equation*}
    x_i=\sqrt{1-\epsilon^2}y_i^{(1)}+ \epsilon y_i^{(2)}
  \end{equation*}
  for all $i\in\{1,\dots,n\}$. Such an ONS can be found by extending $x_1,\dots,x_n$ to an ONS $x_1,\dots,x_n,\tilde x_1,\dots, \tilde x_n$ and defining 
    \begin{equation*}
     y_i^{(1)}=\sqrt{1-\epsilon^2}x_i+ \epsilon\tilde x_i,\ \ y_i^{(2)}=\epsilon x_i -\sqrt{1-\epsilon^2}\tilde x_i.
  \end{equation*}
It is clear that this yields an ONS with
  \begin{equation*}
    \begin{split}
      \sqrt{1-\epsilon^2}y_i^{(1)}+\epsilon y_i^{(2)}&=\sqrt{1-\epsilon^2}(\sqrt{1-\epsilon^2}x_i+\epsilon \tilde x_i)+\epsilon(\epsilon x_i-\sqrt{1-\epsilon^2}\tilde x_i)\\
        &=(1-\epsilon^2)x_i+\epsilon^2x_i=x_i,
      \end{split}
    \end{equation*}
as desired. This establishes the existence of the required ONS.\\
Now define 
  \begin{equation*}
    K=\left\langle y_1^{(2)},y_2^{(2)},\dots, y_n^{(2)}\right\rangle, \ \ \ \  \ \ L=K^\perp.
  \end{equation*}
We claim that there exists an ONS $z_1^{(1)}, z_1^{(2)}, \dots, z_n^{(1)}, z_n^{(2)}$ in $L$ such that 
  \begin{equation*}
    y_i^{(1)}=\tfrac{1}{2}z^{(1)}_i+\tfrac{\sqrt{3}}{2}z^{(2)}_i
    \end{equation*}
    for all $i\in\{1,\dots,n\}$. As above, this is obtained by extending to an ONS $ y_1^{(1)},y_2^{(1)}\dots, y_n^{(1)}, \tilde y_1^{(1)}\dots, \tilde y_n^{(1)}$ in $L$ (possible since $\dim(L)=\infty$) and defining
     \begin{equation*}
     z_i^{(1)}= \tfrac{1}{2}y^{(1)}_i-\tfrac{\sqrt{3}}{2} \tilde y^{(1)}_i,\ \ z_i^{(2)}=\tfrac{\sqrt{3}}{2}y^{(1)}_i+\tfrac{1}{2} \tilde y^{(1)}_i.
  \end{equation*}
Next we claim that there exist isometries $R\in\mathcal{B}(L), V\in\mathcal{B}(K,L)$ such that 
  \begin{equation*}
    R(y_i^{(1)})=z_i^{(1)},\ \ \ V(y_i^{(2)})=z^{(2)}_i, 
    \end{equation*}
and $L=\textup{Im}(V)\oplus\textup{Im}(R)$. To see this, define isometries
  \[
V:K\to L,\quad y^{(2)}_i\mapsto z^{(2)}_i,\qquad
\widetilde R:\langle y^{(1)}_1,\dots,y^{(1)}_n\rangle\to L\ominus V(K),\quad y^{(1)}_i\mapsto z^{(1)}_i.
\]
Since $\dim(L)=\dim(L\ominus V(K))=\infty$, there exists a surjective isometry $R\in\mathcal{B}(L, L\ominus V(K))$ extending $\tilde R$. Thus the desired operator is obtained by considering $R$ as an operator in $\mathcal{B}(L)$.\\
Set
  \begin{equation*}
     I_F=\begin{pmatrix} R & \tfrac{\sqrt{3(1-\epsilon^2)}}{\epsilon}V\\ 0 & id_K \end{pmatrix}.
  \end{equation*}
This is a Brownian unitary and hence a 2-isometry. The proof is completed by showing that $(I_F)_F\to 2id$ in the SOT.\\
For every $x=\sum_{i=1}^n\lambda_ix_i\in F$ (recall $x_i=\sqrt{1-\epsilon^2}y_i^{(1)}+\epsilon y_i^{(2)}$) we have
   \begin{equation*}
     \begin{split}
       I_F(x)&=\sum_{i=1}^n\lambda_i\left(R(\sqrt{1-\epsilon^2}y_i^{(1)})+\tfrac{\sqrt{3(1-\epsilon^2)}}{\epsilon}V(\epsilon y_i^{(2)})+\epsilon id_K(y_i^{(2)})\right)\\
      &=\sum_{i=1}^n\lambda_i\left(\sqrt{1-\epsilon^2}z_i^{(1)}+\sqrt{3(1-\epsilon^2)}z_i^{(2)}+\epsilon y_i^{(2)}\right),
  \end{split}
  \end{equation*}
  and (recall $y^{(1)}_i=\tfrac{1}{2}z^{(1)}_i+\tfrac{\sqrt{3}}{2}z^{(2)}_i$) also
  \begin{equation*}
    \begin{split}
    (2id)(x)&=\sum_{i=1}^n\lambda_i\left(2\sqrt{1-\epsilon^2}y_i^{(1)}+2\epsilon y_i^{(2)}\right)\\
    &=\sum_{i=1}^n\lambda_i\left(\sqrt{1-\epsilon^2} z_i^{(1)}+\sqrt{3(1-\epsilon^2)} z_i^{(2)}+2\epsilon y_i^{(2)}\right) 
      \end{split}
      \end{equation*}
Hence 
  \begin{equation*}
  \|(I_F-2id)(x)\|=\|\epsilon\sum_{i=1}^n\lambda_iy_i^{(2)}\|=\epsilon\|x\|=\tfrac{1}{\dim(F)}\|x\|,
  \end{equation*}
  and therefore
    \begin{equation*}
      \|(I_F-2id)|_F\|\le \tfrac{1}{\dim(F)}.
  \end{equation*}
Now let $G\subset H$ be a finite-dimensional subspace and $\delta>0$. Choose a finite-dimensional subspace $F\subset H$ with $G\subset F$ and $1/\dim(F)<\delta$. Then, for every finite-dimensional subspace $\tilde F$ with $F\subset \tilde F$,
  \begin{equation*}
    \|(I_{\tilde F}-2id)|_G\|\le \|(I_{\tilde F}-2id)|_{\tilde F}\|\le \tfrac{1}{\dim(\tilde F)}\le \tfrac{1}{\dim(F)}<\delta.
  \end{equation*}
Therefore $(I_F)_F\to 2id$ in the SOT. $\hfill\square$
  \ \\
    \begin{Corollary}
  The set of 2-isometries is not SOT-closed.
    \end{Corollary}
  
To prove that the closure is given by all expansive operators, we must address a major obstacle: in general, $T$ does not map orthogonal sets to orthogonal sets. One could try to choose the extension of the ONS $x_1,\dots,x_n$ in the previous proof more carefully, however, there is a convenient workaround. Instead of working solely on $H$, we pass to $H\oplus H$ and $T\oplus T=T^{(2)}$, and then to $H^4$ with $T^{(2)}\oplus T^{(2)}=T^{(4)}$, and construct a net $(I_F)_F$ of $2$-isometries such that $\|(T^{(4)}-I_F )|_G\|\to 0$ for all finite-dimensional subspaces $G$ of the form $\tilde G\oplus 0\oplus 0\oplus 0\subset H^4$. That this is enough to obtain a net of $2$-isometries that converges SOT to $T$ is somehow in the spirit of the Conway and Hadwin result that the SOT-closure with respect to sequences is related to extensions, see \cite{ConHad}.\\
We now record the technical lemmas needed below. It is well known that $2$-isometric operators are expansive, see, for example, \cite[Lemma 2.1]{Ric}. For completeness, we include a proof.

\begin{Lemma}\label{lem:expansive}
Every $2$-isometry $T$ on a Hilbert space is expansive, i.e.\ $\|Tx\|\ge \|x\|$ for all $x$.
\end{Lemma}

\underline{Proof}: \\
Let $T$ be a $2$-isometry. Then
  \[
    T^*T=\tfrac{1}{2}+\tfrac{1}{2}T^*T^*TT\ge\tfrac{1}{2}
  \]
Thus
  \[
    T^*T\ge\tfrac{1}{2}+\tfrac{1}{4}T^*T,
  \]
which implies $T^*T\ge\tfrac{2}{3}$. Recursively we obtain
  \[
    T^*T\ge\tfrac{n}{n+1}
  \]
  for all $n\in\mathbb{N}$, hence $T^*T\ge 1$, which is equivalent to $\|Tx\|\ge\|x\|$ for all $x$.   $\hfill\square$\\

\begin{Lemma}\label{lem:diag-compression}
Let $H$ be a Hilbert space, $T\in\mathcal{B}(H)$ and $F\subset H$ be a finite-dimensional subspace. Then there exists an orthonormal basis $x_1,\dots,x_n$ of $F$ such that $\langle Tx_i,Tx_j\rangle=0$ whenever $i\neq j$.
\end{Lemma}

\underline{Proof}:\\ 
Consider $P_FT^*T|_F$ as a matrix on $F$. Thus, since it is positive, it admits an orthonormal eigenbasis $x_1,\dots,x_n$ of $F$. For $i\neq j$ we have 
  \begin{equation*}
    \langle Tx_i,Tx_j\rangle=\langle T^*Tx_i,x_j\rangle=\langle P_FT^*T|_Fx_i,x_j\rangle=0.\tag*{$\square$}
   \end{equation*}
\vspace{-10pt}
\begin{Lemma}\label{lemOrtho}
Let $H$ be a Hilbert space, $T\in\mathcal{B}(H)$ and $x_1,\dots,x_n$ an ONS such that $Tx_1,\dots,Tx_n$ is an orthogonal family. For $1\ge c\ge0$, define
  \begin{equation*}
    y_i^{(1)}=\sqrt{1-c^2}(x_i\oplus 0)+c(0\oplus x_i), \ \ \ y_i^{(2)}=c(x_i\oplus 0)-\sqrt{1-c^2}(0\oplus x_i)\in H\oplus H,
  \end{equation*}
and set $T^{(2)}=T\oplus T$. Then $\{y_i^{(1)},y_i^{(2)}\}_{i=1}^n$ is an orthonormal set with
 \[
   x_i\oplus 0=\sqrt{1-c^2}\,y_i^{(1)}+c\,y_i^{(2)}\ \textup{and}\ \ \|T^{(2)}y_i^{(1)}\|=\|Tx_i\|.
 \] 
Moreover, 
  \begin{equation*}
    \{T^{(2)}y_i^{(k)};\ k=1,2,\ i=1,\dots,n\}
  \end{equation*}
  forms an orthogonal set.
\end{Lemma}

\underline{Proof}:\\
It is straightforward to verify that $\{y_i^{(1)},y_i^{(2)}\}_{i=1}^n$ is an orthonormal set, $  x_i\oplus 0=\sqrt{1-c^2}\,y_i^{(1)}+c\,y_i^{(2)}$ and $\|T^{(2)}y_i^{(1)}\|=\|Tx_i\|$.\\
To see the remaining orthogonality relations, observe that, for $i\neq j$,
\begin{align*}
\big\langle T^{(2)}y_i^{(1)},T^{(2)}y_j^{(1)}\big\rangle
&=(1-c^2)\langle Tx_i,Tx_j\rangle+c^2\langle Tx_i,Tx_j\rangle=0,\\
\big\langle T^{(2)}y_i^{(1)},T^{(2)}y_j^{(2)}\big\rangle
&=\sqrt{1-c^2}\,c\,\langle Tx_i,Tx_j\rangle
-\sqrt{1-c^2}\,c\,\langle Tx_i,Tx_j\rangle=0,
\end{align*}
and 
\begin{equation*}
\big\langle T^{(2)}y_i^{(1)},T^{(2)}y_i^{(2)}\big\rangle
=\big(\sqrt{1-c^2}\cdot c-c\cdot\sqrt{1-c^2}\big)\,\|Tx_i\|^2=0.\tag*{$\square$}
\end{equation*}

In the proof of Theorem \ref{Theorem1}, we constructed Brownian unitaries. We now replace these with general $2$-isometric operators obtained by dropping the isometric assumption on the upper-right entry in the definition of Brownian unitary.

\begin{Lemma}\label{lem:Brownian-2iso}
Let $R:L\to L$ be an isometry and let $V\in\mathcal{B}(K,L)$ satisfy $R^*V=0$.
Consider $B=\begin{pmatrix}R & V\\ 0 & id_K\end{pmatrix}$ on $L\oplus K$.
Then $B$ is a $2$-isometry.
\end{Lemma}

\underline{Proof}:\\
We have
\[
B^*=\begin{pmatrix}R^*&0\\  V^*&id_K\end{pmatrix},\quad
B^*B=\begin{pmatrix}\mathrm{id}_L&0\\ 0& V^*V+\mathrm{id}_K\end{pmatrix}
\]
and
\[
B^2=\begin{pmatrix}R^2& RV+V\\ 0& id_K\end{pmatrix},\qquad
B^{*2}=\begin{pmatrix}R^{*2}&0\\ V^*R^*+V^* & id_K\end{pmatrix}.
\]
Thus
\begin{align*}
B^{*2}B^2
&=\begin{pmatrix}\mathrm{id}_L& 0\\ 0& V^*R^*RV+V^*RV+ V^*R^*V+V^*V+ \mathrm{id}_K\end{pmatrix}\\
&=\begin{pmatrix}\mathrm{id}_L&0\\ 0& 2V^*V+\mathrm{id}_K\end{pmatrix},
\end{align*}
since $R^*R=\mathrm{id}_L$ and $R^*V=V^*R=0$. Consequently
\[
1-2B^*B+B^{*2}B^2=\begin{pmatrix}0&0\\ 0& -2V^*V+2V^*V+(1-2+1)\,\mathrm{id}_K\end{pmatrix}=\begin{pmatrix} 0 & 0\\ 0 & 0\end{pmatrix},
\]
so $B$ is a $2$-isometry. $\hfill\square$\\
\medskip

  \begin{Theorem}\label{Theorem2}
Let $H$ be an infinite-dimensional Hilbert space. Then
  \begin{equation*}
    \{T\in\mathcal{B}(H);\ \|Tx\|\ge\|x\|\ \forall x\in H\}=\overline{\{T\in\mathcal{B}(H);\ T\ \textup{2-isometry}\}}^{SOT}.
  \end{equation*}
\end{Theorem}
\medskip

\underline{Proof}:\\
The inclusion $\grqq\supset\grqq$ follows from Lemma \ref{lem:expansive} and the observation that if $T_\alpha\to T$ in SOT and $\|T_\alpha x\|\ge\|x\|$ for all $x$ and all $\alpha$, then $\|Tx\|=\lim_\alpha\|T_\alpha x\|\ge\|x\|$.\\
For the converse, let $T\in\mathcal{B}(H)$ be expansive and fix a finite-dimensional subspace $F\subset H$. Set $n=\dim F$ and $\epsilon=1/n$. By Lemma \ref{lem:diag-compression} there exists an orthonormal basis $x_1,\dots,x_n$ of $F$
with $\langle Tx_i,Tx_j\rangle=0$ whenever $i\neq j$.\\
We divide the remainder of the proof into five steps. In the first two steps we apply Lemma \ref{lemOrtho} twice to obtain certain ONSs. In step $3$ we construct a net of $2$-isometries on $H^{(4)}$. In step $4$ we show this net converges $\grqq$partially$\grqq$ to $T^{(4)}$, and in the final step, we use unitary conjugation to obtain the desired net of $2$-isometries on $H$.\\

\smallskip\noindent\emph{Step 1:}\\
On $H^{(2)}=H\oplus H$ with $T^{(2)}=T\oplus T$ define, for $1\le i\le n$,
\[
y_i^{(1)}=\sqrt{1-\epsilon^2}\,(x_i\oplus 0)\ +\ \epsilon\, (0\oplus x_i),\qquad
y_i^{(2)}=\epsilon\,(x_i\oplus 0)\ -\ \sqrt{1-\epsilon^2}\, (0\oplus x_i).
\]
By Lemma \ref{lemOrtho}, $\{y_i^{(1)},y_i^{(2)}\}_{i=1}^n$ is an orthonormal set, 
  \[
    x_i\oplus 0=\sqrt{1-\epsilon^2}\,y_i^{(1)}+\epsilon\,y_i^{(2)}, 
  \]    
and
  \[
    \|T^{(2)}y_i^{(1)}\|=\|Tx_i\|.
  \]
Moreover,
  \begin{equation*}
    \{T^{(2)}y_i^{(k)};\ k=1,2,\ i=1,\dots,n\}
  \end{equation*}
is an orthogonal set.

\smallskip\noindent\emph{Step 2:}\\
Apply Lemma \ref{lemOrtho} once more to $T^{(2)}$ and $y_1^{(1)},\dots,y_n^{(1)}$. Define
\[
\hat y_i^{(1)}=y_i^{(1)}\oplus 0\oplus 0,\qquad \hat y_{i+n}^{(1)}=0\oplus 0\oplus y_i^{(1)},\qquad \hat y_i^{(2)}=y_i^{(2)}\oplus 0\oplus 0,
\]
and, note that $\tfrac{1}{\|Tx_i\|}\le 1$ since $T$ is expansive,
\[
z_i^{(1)}=\frac{1}{\|Tx_i\|}\,\hat y_i^{(1)}+\sqrt{1-\frac{1}{\|Tx_i\|^2}}\,\,\hat y_{i+n}^{(1)},\quad
z_i^{(2)}=\sqrt{1-\frac{1}{\|Tx_i\|^2}}\,\,\hat y_i^{(1)}-\frac{1}{\|Tx_i\|}\,\hat y_{i+n}^{(1)}.
\]
Then $\{z_i^{(1)},z_i^{(2)}\}_{i=1}^n$ is an orthonormal set,
\[
\hat y_i^{(1)}=\tfrac{1}{\|Tx_i\|}z_i^{(1)}+\sqrt{1-\tfrac{1}{\|Tx_i\|^2}}z_i^{(2)},
\]
\[
\|Tx_i\|=\|T^{(2)}y_i^{(1)}\|=\|T^{(4)}z_i^{(1)}\|,
\]
and
  \begin{equation}
    \{T^{(4)}z_i^{(k)},\ 1\le k\le 2,\ 1\le i\le n\} 
    \label{EqOrtho}
  \end{equation}
is an orthogonal set.\\

\smallskip\noindent\emph{Step 3:}\\
Set
\[
K=\operatorname{span}\{\hat y^{(2)}_1,\dots,\hat y^{(2)}_n\}\subset H^{(4)},\qquad L=K^\perp.
\]
It holds that, $\text{for all }i,j\in\{1,\dots,n\}$,
  \begin{equation*}
    \begin{split}
      \langle T^{(4)}z_i^{(1)},\,\hat y^{(2)}_j\rangle\;&=\tfrac{1}{\|Tx_i\|}\langle T^{(2)}y_i^{(1)},y_j^{(2)}\rangle\\
       &=\tfrac{1}{\|Tx_i\|}(\langle (\sqrt{1-\epsilon^2}Tx_i)\oplus (\epsilon Tx_i), (\epsilon x_j)\oplus (-\sqrt{1-\epsilon^2}x_j)\rangle\\
       &=\tfrac{1}{\|Tx_i\|}\sqrt{1-\epsilon^2}\epsilon(\langle Tx_i,x_j\rangle-\langle Tx_i,x_j\rangle)=\;0\\
    \end{split}
  \end{equation*}
and similarly
  \[ 
     \langle T^{(4)}z_i^{(2)},\,\hat y^{(2)}_j\rangle\;=\sqrt{1-\tfrac{1}{\|Tx_i\|}^2}\langle T^{(2)}y_i^{(1)},y_j^{(2)}\rangle\\
    =\;0,
  \]
hence $T^{(4)}z_i^{(k)}\in L$ for $k=1,2$. For
  \[
    \sigma_i=\tfrac{1}{\epsilon}\sqrt{(1-\epsilon^2)(1-\tfrac{1}{\|Tx_i\|^2})}, 
  \]  
define an operator $V:K\to L$ by 
\[
V(\hat y^{(2)}_i)=\sigma_iT^{(4)}z^{(2)}_i.
\]
By (\ref{EqOrtho}), we can define the isometric operator
\[
R_0:\ \operatorname{span}\{\hat y^{(1)}_1,\dots,\hat y^{(1)}_n\}\ \longrightarrow L\ominus V(K),
\qquad
R_0(\hat y^{(1)}_i)=\frac{T^{(4)}z^{(1)}_i}{\|Tx_i\|}.
\]
Since $\dim(V(K))<\infty$, we have that $\dim(L)=\dim(L\ominus V(K))=\infty$, and since $\hat y^{(1)}_i\in L$, $R_0$ extends to an isometry
$R:L\to L\ominus V(K)$.\\
View $R$ as an operator in $\mathcal{B}(L)$ and set
\begin{equation}\label{eq:IF-def}
I_F\ =\ \begin{pmatrix} R &  V\\[2pt] 0 & \mathrm{id}_K\end{pmatrix}
\quad\text{on }L\oplus K=H^4,
\end{equation}
By Lemma \ref{lem:Brownian-2iso}, $I_F$ is a $2$-isometry.\\

\smallskip\noindent\emph{Step 4:}\\
Let $x\in F$ and write $x=\sum_{i=1}^n\lambda_i x_i$. Identify $x_i$ with $x_i\oplus0\oplus0\oplus0$. Recall that $x_i\oplus 0=\sqrt{1-\epsilon^2}\,y_i^{(1)}+\epsilon\,y_i^{(2)}$. Using the definition of $R$ and $V$,
  \begin{align*}
    I_F(x_i)
    &=\sqrt{1-\epsilon^2}I_F(\hat y_i^{(1)})+\epsilon I_F(\hat y_i^{(2)})\\
    &=\sqrt{1-\epsilon^2}R(\hat y_i^{(1)})+\epsilon V(\hat y_i^{(2)})+ \epsilon \hat y_i^{(2)}\\
    &=\sqrt{1-\epsilon^2}\frac{T^{(4)}z^{(1)}_i}{\|Tx_i\|}+\epsilon\sigma_iT^{(4)}z^{(2)}_i+\epsilon\hat y_i^{(2)}\\
    &=\sqrt{1-\epsilon^2}\Big(\frac{1}{\|Tx_i\|}T^{(4)}z^{(1)}_i+\sqrt{1-\frac{1}{\|Tx_i\|^2}}\,T^{(4)}z_i^{(2)}\Big)+\epsilon\hat y_i^{(2)}.
\end{align*}
On the other hand, recall that $\hat y_i^{(1)}=\tfrac{1}{\|Tx_i\|}z_i^{(1)}+\sqrt{1-\tfrac{1}{\|Tx_i\|^2}}z_i^{(2)}$,
  \begin{align*}
    T^{(4)}(x_i)
    &=\sqrt{1-\epsilon^2}\,T^{(4)}\hat y_i^{(1)}+\epsilon\,T^{(4)}\hat y_i^{(2)}\\
    &=\sqrt{1-\epsilon^2}\Big(\frac{1}{\|Tx_i\|}T^{(4)}z_i^{(1)}+\sqrt{1-\frac{1}{\|Tx_i\|^2}}T^{(4)}z_i^{(2)}\Big)+\epsilon\,T^{(4)}\hat y_i^{(2)}.
\end{align*}
Subtracting gives
  \[
    \big(T^{(4)}-I_F\big)x_i= \epsilon\Big(T^{(4)}-id\Big)\hat y_i^{(2)}.
  \]
Using the orthogonality of $\hat y_1^{(2)}\dots,\hat y_n^{(2)}$, we obtain
  \[
    \|(T^{(4)}-I_F)x\|=\epsilon\|(T^{(4)}-id)(\sum_{i=1}^n\lambda_i\hat y_i^{(2)})\|\le \epsilon(\|T\|+1)\|x\|.
  \]
Hence
\begin{equation*}
\|(T^{(4)}-I_F)|_{F\oplus0\oplus0\oplus0}\|\ \le\ \epsilon\,(\|T\|+1)\ =\ \frac{\|T\|+1}{\dim(F)}.
\end{equation*}
\smallskip\noindent\emph{Step 5:}\\
Since $H$ is infinite-dimensional, there exists a unitary $U_F:H\to H^{(4)}$ with $U_F x=x\oplus0\oplus0\oplus0$ for all $x\in \textup{span}(F\cup T(F))$. Set $J_F=U_F^*I_FU_F$. Then $J_F$ is a $2$-isometry on $H$, and
  \[
    T^{(4)}(U_F(x))=(Tx)\oplus0\oplus 0\oplus 0=U_F(Tx)
  \]
for all $x\in F$, so that
  \[
    U_F^*T^{(4)}U_F|_F=T|_F.
  \]
It remains to show that the net $(J_F)_F$ converges SOT to $T$. For this, let $\delta>0$ and $G\subset H$ be a finite-dimensional subspace. Choose a finite-dimensional subspace $F\subset H$ such that $G\subset F$ and $1/\dim(F)\le \delta$. Then, for all finite-dimensional subspaces $\tilde F\subset H$ containing $F$,
  \[
    \begin{split}
      \|(T-J_{\tilde F})|_{G}\|\le \|(T-J_{\tilde F})|_{\tilde F}\|&=\|U_{\tilde F}^*(T^{(4)}-I_{\tilde F})U_{\tilde F}|_{\tilde F}\|\\
      &=\|(T^{(4)}-I_{\tilde F})|_{\tilde F\oplus0\oplus0\oplus0}\|\\\
      & \le\ \frac{\|T\|+1}{\dim(\tilde F)}\\
       & \le\ \frac{\|T\|+1}{\dim(F)}\le \delta(\|T\|+1).
    \end{split}
  \]
Hence $J_F\to T$ in the SOT.
\hfill$\square$\\

\begin{Corollary}
Let $H$ be an infinite-dimensional Hilbert space. Then the set of $2$-isometric operators is WOT-dense in $\mathcal{B}(H)$.
\end{Corollary}

\underline{Proof}:\\
It is well known that the WOT-closure of the set of unitary operators equals all contractive operators. Hence an arbitrary operator $T$ is the WOT-limit of a net of the form $(t U_\alpha)_\alpha$, where $1\le t\in\mathbb{R}$ and $(U_\alpha)_\alpha$ are unitary operators converging WOT to $\tfrac{T}{t}$. By Theorem \ref{Theorem2}, each $t U_\alpha$ lies in the SOT-closure of the set of $2$-isometric operators, finishing the proof. $\hfill\square$\\

\begin{Remark}\label{Remark}
It is natural to ask about the SOT-closure of the set of $3$-isometries, that is 
  \[
    \{T\in\mathcal{B}(H);\ -1+3T^*T-3(T^*)^2T^2+(T^*)^3T^3=0\}.
  \]
  It turns out that for every $2$-nilpotent operator $A$, the operator $id+A$ is a $3$-isometry, see \cite{Rus} or \cite{BMN}. Using the aforementioned result by Halmos that the set of $2$-nilpotent operators is SOT-dense in $\mathcal{B}(H)$, it follows that the SOT-closure of the set of $3$-isometries is $\mathcal{B}(H)$.
  \end{Remark}

\section*{Acknowledgement}
A big thank you goes to Pawel Pietrzycki for pointing out this problem to me. I also wish to thank Michael Hartz for pointing out Remark \ref{Remark}.
\ \\
\ \\

\bibliography{bib_SOTAbschluss} 
\bibliographystyle{plain}

\textit{Email address:} scherer@math.uni-sb.de

\end{document}